\documentclass[twocolumn]{autart_edited}    

\usepackage{graphicx} 
\usepackage{amsmath}
\usepackage{amssymb}
\usepackage{color}

\newcommand{\AAA}{{\mathcal A}}

\newcommand{\AC}{{\mathcal A}^{\mathsf C}}

\newcommand{\II}{{\mathbb I}}

\newcommand{\DA}{\partial \mathcal{A}}

\newcommand{\DAM}{\left [\DA\right]_{\mathcal{-}}}

\newcommand{\ee}{\varepsilon}

\newcommand{\cl}{\mathsf{cl}}

\newcommand{\Int}{\mathsf{int}}

\newcommand{\NN}{{\mathbb N}}

\newcommand{\RR}{{\mathbb R}}

\newcommand{\UU}{{\mathcal U}}

\newcommand{\qq}{{\mathbf q}}
\newcommand{\hh}{{\mathbf h}}
\newcommand{\tth}{{\boldsymbol \theta}}
\newcommand{\ds}{\displaystyle}

\newtheorem{definition}{Definition}
\newtheorem{proposition}{Proposition}
\newtheorem{remark}{Remark}

\begin{document}

\begin{frontmatter}
\runtitle{Barrier, Safe Set and Pendulum with Non-rigid Cable}  
                                             
\title{Barriers and Potentially Safe Sets in Hybrid Systems: Pendulum with Non-Rigid Cable\thanksref{footnoteinfo}} 

\thanks[footnoteinfo]{© 2016. This manuscript version is made available under the CC-BY-NC-ND 4.0 license http://creativecommons.org/licenses/by-nc-nd/4.0/. DOI: 10.1016/j.automatica.2016.07.001. Corresponding author: Jean L\'{e}vine.}

\author[CAS]{Willem Esterhuizen}\ead{willem.esterhuizen@mines-paristech.fr},    
\author[CAS]{Jean L\'{e}vine}\ead{jean.levine@mines-paristech.fr},             
\address[CAS]{CAS, Systems and Control Centre, MINES ParisTech, PSL Research University, 60 Bd Saint-Michel, 75006 Paris, France.}   
          

\begin{abstract}                       
This paper deals with an application of the notion of barrier in mixed constrained nonlinear systems to an example of a pendulum mounted on a cart with non-rigid cable, whose dynamics may switch to free-fall when the tension of the cable vanishes. We present a direct construction of the boundary of the potentially safe set in which there always exists a control such that the cable never goes slack. A discussion on the dependence of this set with respect to the pendulum and cart masses is then sketched.
\end{abstract}

\end{frontmatter}

\section{Introduction}
This paper presents an application and slight extension of the recent work on \emph{barriers} in constrained nonlinear systems, see \cite{DeDona_siam,EsterhuizenPHDThesis,Ester_Lev_arxiv}. Given a pendulum on a cart with the rigid bar replaced by a massless cable, we aim at designing a control law which guarantees that the cable always remains taut. The study of this system may be useful to the investigation of safely controlling overhead cranes where slackness of the cable would result in free-fall of the working mass, which would therefore be uncontrolled, and thus potentially harmful for the system and its environment. Such a system whose dynamics may switch conditionally to an event which is, itself, a function of the state and input, is generally called a \emph{hybrid system} (see e.g. \cite{Mitchell2003,VanDerShaft2000,Lygeros2007}). The reader may also refer to \cite{KissPHDThesis,KLM-pp,KLM-jss} for studies on modeling and trajectory planning of \emph{weight handling equipment}. A similar problem appears in \cite{NicotraNalGor_IFACE2014} where the authors study \emph{tethered unmanned aerial vehicles} in the different perspective of designing a stabilizing feedback controller. 

For a constrained nonlinear control system, the \emph{admissible set} is the set of all initial conditions for which there exists a control such that the constraints are satisfied for all time. Under mild assumptions, this set is closed and its boundary consists of two complementary parts. One of them, called the \emph{barrier}, enjoys the so-called \emph{semi-permeability} property \cite{Isaacs} and its construction is done via a minimum-like principle \cite{DeDona_siam,EsterhuizenPHDThesis,Ester_Lev_arxiv}. Our approach to solving the above mentioned problem of the pendulum on a cart is to find this system's admissible set and to guarantee the cable tautness as follows: if the state remains in the admissible set's interior, the control can be arbitrary in some state-dependent constraint set for almost all time and, if the state reaches the barrier, a special control, which we indeed exhibit, needs to be employed in order to keep the cable taut. This admissible set may be interpreted as a \emph{safe set}, or more precisely as a \emph{potentially safe set}. 

Note also that we emphasize on systems with \emph{mixed} constraints, i.e. constraints that are functions, in a coupled way, of the control and the state \cite{Clarke_DiPinho,Hestenes,Ester_Lev_arxiv}, the reason being that tautness of the cable, which is expressed by the fact that the tension in the cable remains nonnegative, can be shown to be equivalent to imposing a mixed constraint. Such constraints are by far more complicated than pure state constraints since they are control dependent, with controls that may be discontinuous with respect to time, thus possibly creating jumps on the constraint set.

Admissible sets are strongly related to invariant sets \cite{Chutinan03c,Teel:HybridBook2012} and viability kernels \cite{Aubin,Lygeros99,Tomlin2000,VanDerShaft2000,Kaynama:2012:CVK:2185632.2185644,Mitchell2003,Mitchel2005,Lhommeau:capture:11}. Our approach contrasts with these works by the fact that in place of computing flows or Lyapunov functions or solutions of Hamilton-Jacobi equations over the whole domain, \emph{we reduce the computations to the boundary of the set under study}. The same kind of comparison also holds with barrier Lyapunov functions \cite{Tee2009a}, or barrier certificates \cite{Prajna06}.

The originality of the results of this paper is threefold: 
\begin{itemize}
\item the interpretation of the cable tautness / slackness as a mixed constraint may be found in \cite{NicotraNalGor_IFACE2014} but, as already said, with a different stabilization objective.
In this paper, we are interested in the analysis and computation of the associated admissible set, namely the largest state domain where one can find an open-loop control such that the cable remains taut, which is new to the authors' knowledge; 
\item the computation of this admissible set by focusing on its boundary strongly contrasts, in spirit, with the various theoretical constructions found in the literature \cite{Aubin,Lygeros99,Tomlin2000,VanDerShaft2000,Kaynama:2012:CVK:2185632.2185644,Mitchell2003,Mitchel2005,Lhommeau:capture:11,NicotraNalGor_IFACE2014} where numerical integration is used to compute flows, each step being simple but the number of steps and iterations exponentially increasing with the dimension of the problem;
\item the necessary conditions used here have been obtained in \cite{Ester_Lev_arxiv} at the exception of the  terminal condition called \emph{ultimate tangentiality condition}. This new terminal condition, introduced to overcome a double problem of singularity and nonsmoothness, is essential for the computation of the barrier: the latter equations cannot be integrated without suitable terminal conditions and the ultimate tangentiality condition of \cite{Ester_Lev_arxiv} turned out to be too coarse to obtain a solution.
\end{itemize}

The paper is organised as follows. In Section \ref{sec:Barriers_In_Constrained_Sys} we summarise the main results from \cite{DeDona_siam}, \cite{EsterhuizenPHDThesis} and \cite{Ester_Lev_arxiv} which we present without proofs.  In Section \ref{sec:Barrier_for_Pend} we construct the system's barrier. Section \ref{sec:Discussion} provides a discussion of the physical interpretations of the results, and the paper ends with Section \ref{sec:Conclusions} that summarises the conclusions and points out future research.

\section{Barriers in Nonlinear Control System with Mixed Constraints}\label{sec:Barriers_In_Constrained_Sys}
\subsection{Constrained Nonlinear Systems with Mixed Constraints}
The contents of this section is borrowed from \cite{EsterhuizenPHDThesis} and \cite{Ester_Lev_arxiv}, where more details may be found. However, Proposition~\ref{ult-tan-1d-pr} and Theorem~\ref{BarrierTheorem1} of this paper slightly extend the ones of these references. We consider the following nonlinear system with mixed constraints:
\begin{align}
	\label{eq:state_space}
	& \dot{x} =  f(x,u), \\
	\label{eq:initial_condition}
	& x(t_0) = x_0, \\
	\label{eq:input_constraint}
	& u  \in \UU, \\
	\label{eq:state_const}
	& g_i\big(x(t), u(t)\big)  \leq   0 \quad  \mathit{a.e.~} t \in [t_0, \infty) \quad i=1,...,p
\end{align}
where $x(t)\in \RR^{n}$. 
The set $\UU$ is the set of Lebesgue measurable functions from $[t_0, \infty)$ to  $U$, a given compact convex subset of $
^{m}$; Thus $u$ is a measurable function such that $ u(t) \in U$ for almost all $t\in [t_0, \infty)$.  

We denote by $x^{(u,x_0,t_0)}(t)$ the solution of the differential equation~\eqref{eq:state_space} at $t$ with input \eqref{eq:input_constraint}  and initial condition~\eqref{eq:initial_condition}. Sometimes the initial time or initial condition need not be specified, in which cases we will use the notation $x^{(u,x_0)}(t)$ or $x^u(t)$ respectively.

The constraints \eqref{eq:state_const}, called \emph{mixed constraints} \cite{Clarke_DiPinho,Hestenes}, explicitly depend both on the state and the control. We denote by $g(x,u)$ the vector-valued function whose $i$-th component is $g_i(x,u)$. By $g(x,u)\prec 0$ (resp. $g(x,u)\preceq 0 $) we mean $g_i(x,u) < 0$ (resp. $g_i(x,u) \leq 0$) for all $i$. By $g(x,u)\circeq 0$, we mean $g_i(x,u) = 0$ for at least one $i$. As said before, even if $g$ is smooth, the mapping $t\mapsto g(x(t),u(t))$ is only measurable and the associated mixed constraints are thus assumed to be  satisfied almost everywhere.

\subsection{The Admissible Set}
We define the following sets:
\begin{gather} \ds 
	G \triangleq \bigcup_{u\in U} \{x\in\RR^n : g(x,u)\preceq 0\} \label{def:G}
	\\ \ds 
	G_0 \triangleq \{x \in G : \min_{u\in U} \max_{i\in\{1,...,p\}}g_i(x,u)=0 \}\label{def:G0} 
	\\ \ds 
	G_{-} \triangleq \bigcup_{u\in U} \{x\in\RR^n : g(x,u) \prec 0\} \label{def:G_-}
\end{gather}

We further assume:
\begin{description}
	\item[(A2.1)] $f$ is an at least $C^{2}$ vector field of $\RR^{n}$ for every $u$ in an open subset $U_1$ of $\Rset^{m}$ containing $U$, whose dependence with respect to $u$ is also at least $C^{2}$.
	\item[(A2.2)] There exists a constant $0 < C < +\infty$ such that the following inequality holds true:
	$$\sup_{u\in U}\vert x^T f(x,u) \vert \leq C(1+ \Vert x \Vert^{2} ), \quad \mbox{\textrm{for all}}~x$$
	where the notation $x^T f(x,u)$ indicates the inner product of the two vectors $x$ and $f(x,u)$.
	\item[(A2.3)] The set $f(x,U)$, called the \emph{vectogram} in \cite{Isaacs}, is convex for all $x\in \RR^{n}$.
	\item[(A2.4)] $g$ is at least $C^{2}$ from $\RR^{n}\times U_1$ to $\RR^p$ and convex with respect to $u$ for all $x\in \RR^{n}$. 
 \end{description}

We also introduce the following state-dependent control set:
\begin{equation}
U(x) \triangleq \{u \in U : g(x,u) \preceq 0 \}\quad \forall x \in G.\label{def:U(x)}
\end{equation}
The convexity of $U$ and (A2.4) imply that $U(x)$ is convex for all $x\in G$ and, since $g$ is continuous, the multivalued mapping $x\mapsto U(x)$ is closed with range in the compact set $U$, and therefore upper semi-continuous (u.s.c.) (see e.g. \cite{berge,F}).

We assume that, for every $x\in G$, the set $U(x)$ is locally expressible as 
\begin{equation}\label{def:UStraight}
U(x) \triangleq \{u\in\RR^m : \gamma_i(x,u)\leq 0, i = 1,\dots,r\}
\end{equation}
the functions $\gamma_{i}$ being of class $C^2$, linearly independent, and convex with respect to $u$ for all $x\in G$. 

For a pair $(x,u)\in \RR^n \times U$, we denote by $\II(x,u)$ the set of indices, possibly empty, corresponding to the ``active'' mixed constraints:
\begin{equation}\label{IIdef}
\II(x,u) \triangleq \{ i\in \{1,\ldots,r\} : \gamma_{i}(x,u) = 0\}.
\end{equation}
The number $\#(\II(x,u))$ of elements of $\II(x,u)$ thus represents the number of ``active'' constraints among the $r$ independent constraints at $(x,u)$. We further assume:

\begin{description}
\item[(A2.5)] For almost all $z$ in a neighborhood of $G_0$ and all $u\in U(z)$ such that $0<\#(\II(z,u))$,  the (row) vectors 
$\frac{\partial\gamma_{i}}{\partial u}(z,u)$, $i \in \II(z,u)$, are linearly independent.
\end{description}

\begin{definition}[Admissible Set]\emph{\cite{Ester_Lev_arxiv}}
	\label{def:admiss_states}
	We say that a point $x_0$ is \emph{admissible} if, and only if, there exists, at least, one input function $u\in \UU$, such	that~\eqref{eq:state_space}--\eqref{eq:state_const} are
	satisfied:
	\begin{equation}\label{eq:Admiss_states}
	\AAA \triangleq \{x_0 \in G: \exists u\in \UU,~ g\big(x^{(u,x_{0})}(t), u(t)\big) \preceq 0, \mathit{a.e.~} t\}.
	\end{equation}
\end{definition}

As in \cite{PBGM} a \emph{Lebesgue} point, for a given control $u\in \UU$ is a time $t\in [t_0,\infty)$ where $u$ is continuous, the interval $[t_0,\infty)$ being possibly deprived of a bounded subset of zero Lebesgue measure which does not contain $t$.

If $u_{1}\in \UU$ and $u_{2}\in \UU$, and if $\tau\geq t_{0}$ is given, the concatenated input $v$, defined by $v(t)= \left\{ \begin{array}{ll} u_{1}(t)&\mbox{\textrm if~} t\in [t_{0}, \tau[\\u_{2}(t)&\mbox{\textrm if~} t \geq \tau\end{array}\right.$ satisfies $v\in \UU$. The concatenation operator relative to $\tau$ is denoted by $\Join_{\tau}$, i.e. $v=u_{1}\Join_{\tau} u_{2}$.

Since system (\ref{eq:state_space}) is time-invariant, the initial time $t_0$ may be taken as 0. When clear from the context, ``$\forall t$'' or ``for \emph{a.e} $t$'' will mean
``$\forall t \in [0, \infty)$'' or ``for \emph{a.e.} $t\in [0, \infty)$'', where \emph{a.e.} is understood with respect to the Lebesgue measure.

\begin{proposition}\label{close-cor}\emph{\cite{Ester_Lev_arxiv}}
	Under assumptions (A2.1) - (A2.4) the set $\AAA$ is closed.
\end{proposition}

\begin{remark}\label{rem:A22}
Assumption (A2.2), which implies an at most linear growth of $f$ with respect to $x$, is introduced to guarantee the uniform boundedness  and uniform convergence of a sequence of integral curves that appear in the proof of Proposition~\ref{close-cor}. However, this condition is far from being necessary and many systems
do not satisfy it though having bounded trajectories under admissible controls. Any other condition on $f$ ensuring uniform boundedness (see e.g. \cite{F}) would give similar compactness results.
\\
Assumption (A2.5) is used in the proof of Theorem \ref{BarrierTheorem1} (see Subsection \ref{barrierTh:subsec}). It replaces the stringent independence condition (A4) of \cite{Ester_Lev_arxiv} which is not satisfied in many examples, including the pendulum one of Section \ref{sec:Barrier_for_Pend}.
\end{remark}

We denote by $\DA$ the boundary of $\AAA$ and $\AC$ its complement. We indeed have $\partial \AAA \subset \AAA$.
\begin{definition}\label{def:barrier}
	The set $\DAM = \DA\cap G_-$ is called the \emph{barrier} of the set $\AAA$.
\end{definition}
It is characterised by the two next Propositions, proved in \cite{Ester_Lev_arxiv}, and the (new) Proposition~\ref{ult-tan-1d-pr} of Subsection~\ref{ult-tan:subsec}. 
\begin{proposition}\label{boundary:prop}\emph{\cite{Ester_Lev_arxiv}}
	Assume (A2.1) to (A2.4) hold. The barrier $\DAM$ is made of points $\bar{x}\in G_-$ for which there exists $\bar{u}\in\UU$ and an integral curve $x^{(\bar{u},\bar{x})}$ entirely contained in $\DAM$ either until it intersects $G_0$, i.e. at a point $z = x^{(\bar{u},\bar{x})}(\bar{t})$, for some $\bar{t}$, such that $\min_{u\in U} \max_{i=1,\ldots,p} g_{i}(z,u) = 0$, or that never intersects $G_0$.
\end{proposition}

\begin{proposition}\label{bar-sem-cor}\emph{\cite{Ester_Lev_arxiv}}
	Assume (A2.1) to (A2.4) hold. Then from any point on the boundary $\DAM$, there cannot exist a trajectory penetrating the interior of $\AAA$ before leaving $G_{-}$.
\end{proposition}

\subsection{Barrier End Point Condition}\label{ult-tan:subsec}

We introduce the notation 
\begin{equation}\label{gt:def}
\tilde{g}(x) \triangleq  \min_{u\in U} \max_{i = 1,\dots, p} g_i(x,u),
\end{equation}
i.e. $G_0 = \{ x : \tilde{g}(x) = 0 \}$. In \cite{EsterhuizenPHDThesis} and \cite{Ester_Lev_arxiv} it was shown that $\tilde{g}$ is locally Lipschitz (this is a version of Danskin's Theorem, see e.g. \cite{danskin}), and therefore differentiable almost everywhere. 

The intersection between $\DAM$ and $G_0$, if it exists, must satisfy the condition given in the next proposition.
\begin{proposition}\label{ult-tan-1d-pr}
	Assume (A2.1) to (A2.4) hold. Consider $\bar{x} \in \DAM$ and $\bar{u}\in \UU$ as in Proposition~\ref{boundary:prop}, i.e. such that the integral curve $x^{(\bar{u},\bar{x})}(t) \in \DAM$ for all $t$ in some time interval until it reaches $G_{0}$ at some finite time $\bar{t}\geq 0$. Then, the point $z= x^{(\bar{u},\bar{x})}(\bar{t})\in \cl(\DAM)\cap G_{0}$,  satisfies
\begin{equation}\label{ineq:ult_tan}
		\min_{u\in U(x^{(\bar{u},\bar{x})}(\bar{t}_-))} D\tilde{g}(x^{(\bar{u},\bar{x})}(\bar{t}_-)) f(x^{(\bar{u},\bar{x})}(\bar{t}), u) \geq 0,
	\end{equation}
	where $D\tilde{g}$ is the gradient of $\tilde{g}$, $h(x(\bar{t}_-))$ indicating the left limit of $h(x(\tau))$, when $\tau \nearrow \bar{t}$ (i.e. with $\tau < \bar{t}$), of an arbitrary function or multivalued mapping $h$, not necessarily continuous.
	
	Moreover, if the point $z$ is a differentiability point of $\tilde{g}$, condition \eqref{ineq:ult_tan} reads
\begin{equation}\label{eq:ult_tan_smooth}
\min_{u\in U(z)} D\tilde{g}(z) f(z, u) = 0
\end{equation}

\end{proposition}
\begin{proof} 
	\sloppy Let $x_0\in \DAM$, then there exists a control $\bar{u}\in \UU$ such that $\tilde{g}(x^{(\bar{u},x_0)}(t)) < 0$ until $x^{(\bar{u},x_0)}$ intersects $G_0$ at some $\tilde{t}$ that we assume finite. Consider an open set ${\mathcal O} \subset \RR^n$ such that $x_0 +\ee h \in \AC$, the complement of $\AAA$, for all $h \in {\mathcal O}$ and $\Vert h\Vert \leq H$, with $H$ arbitrarily small, and all $\ee$ sufficiently small.
	
	Introduce a needle perturbation of $\bar{u}$, labelled $u_{\kappa,\ee}$ as in Appendix \ref{Appendix:var}, 
	where $v \in U(x^{(\bar{u},x_0)}(\tau))$ for all $t\in [\tau-l\ee, \tau[$, at some Lebesgue point $\tau$ of $\bar{u}$ before $x^{(\bar{u},x_0)}$ intersects $G_0$.
Because $x_0 + \ee h \in \AC$, $\exists t_{\ee,\kappa,h}<\infty$ at which  $x^{(u_{\kappa,\ee},x_0+\ee h)}(t_{\ee,\kappa,h})$ crosses $G_0$, see Proposition \ref{bar-sem-cor}. As a result of the uniform convergence of $x^{(u_{\kappa,\ee},x_0+\ee h)}$ to $x^{(\bar{u},x_0)}$, there exists a $\bar{t} \geq \tilde{t}$, s.t. $x^{(u_{\kappa,\ee},x_0+\ee h)}(t_{\ee,\kappa,h})\rightarrow x^{(\bar{u},x_0)}(\bar{t})$ as $\ee \rightarrow 0$ and, according to the continuity of $\tilde{g}$, we have
	$$
	\lim_{\ee \rightarrow 0} \tilde{g}(x^{(u_{\kappa,\ee},x_0+\ee h)}(t_{\ee,\kappa,h})) = 0 = \tilde{g}(x^{(\bar{u},x_0)}(\bar{t})).
	$$
	Because $\tilde{g}(x^{(u_{\kappa,\ee},x_0+\ee h)}(t_{\ee,\kappa,h})) = 0$ and $\tilde{g}(x^{(\bar{u},x_0)}(t_{\ee,\kappa,h})) \leq 0$ (recall that $\tilde{g}(x^{(\bar{u},x_0)}(t_{\ee,\kappa,h})) \leq g(x^{(\bar{u},x_0)}(t_{\ee,\kappa,h}), \bar{u}(t_{\ee,\kappa,h})) \leq 0$ since the pair $(x^{(\bar{u},x_0)}(t),\bar{u}(t))$ satisfies the constraints for all $t$), we have that
	$$
	\tilde{g}(x^{(u_{\kappa,\ee},x_0+\ee h)}(t_{\ee,\kappa,h})) - \tilde{g}(x^{(\bar{u},x_0)}(t_{\ee,\kappa,h})) \geq 0.
	$$
	
	Recall from Appendix \ref{Appendix:var} that	
	\begin{equation}
		\begin{aligned}
		x^{(u_{\kappa,\ee},x_0+\ee h)}(t_{\ee,\kappa,h}) = &x^{(\bar{u},x_0)}(t_{\ee,\kappa,h}) \\ &\,\,+ \ee w(t_{\ee,\kappa,h},\kappa,h) + O(\ee^2)
		\end{aligned}
	\end{equation}
	where $w(t,\kappa,h)$ satisfies \eqref{needle-eq}.
	
	Since $\tilde{g}$ is almost everywhere differentiable, we have:
	\begin{equation}\label{genDerOfGTilde}
		\begin{aligned}
		&\frac{\tilde{g}(x^{(u_{\kappa,\ee},x_0+\ee h)}(t_{\ee,\kappa,h})) - \tilde{g}(x^{(\bar{u},x_0)}(t_{\ee,\kappa,h}))}{\ee} \\ &=D \tilde{g}(x^{(\bar{u},x_0)}(t_{\ee,\kappa,h})) . w(t_{\ee,\kappa,h},\kappa,h)) +O(\ee) \geq 0
		\end{aligned}
	\end{equation}
	for every $v \in U(x^{(\bar{u},x_0)}(\tau))$ and almost every $\ee$ and $h$.
	
	Note that $\tilde{g}(x^{(\bar{u},x_0)}(t)) \prec 0$ for all $t\in ]\bar{t}-\eta,\bar{t}[$, which implies, according to \eqref{gt:def} and (A2.4), that there exists open sets $V(t)$ such that $\cl(V(t)) \subset U(x^{(\bar{u},x_0)}(t))$ and $g(x^{(\bar{u},x_0)}(t),\bar{u}(t)) \preceq g(x^{(\bar{u},x_0)}(t),v) \prec 0$ for all $v\in V(t)$ and all $t\in ]\bar{t}-\eta,\bar{t}[$. Moreover, the multivalued mapping $t\mapsto \cl(V(t))$ may be chosen lower semi-continuous on $]\bar{t}-\eta,\bar{t}[$ (see e.g. \cite{Repov_etal_98} and the survey \cite{Repov_arxiv_14}).
Therefore we can select a continuous selection $v_{\tau} \in \cl(V(\tau)) \subset U(x^{({\bar{u},\bar{x}})}(\tau))$ for $\tau \in ]\bar{t} - \eta, \bar{t}[$ (see again \cite{Repov_arxiv_14}) such that $\lim_{\tau\nearrow \bar{t}} v_{\tau} = v$.  Taking the limit as $\tau \nearrow \bar{t}$ \eqref{genDerOfGTilde} becomes:
	\[
	\min_{v\in U(x^{(\bar{u},\bar{x})}(\bar{t}_-))} D\tilde{g}(x^{(\bar{u},\bar{x})}(\bar{t}_-)) f(x^{(\bar{u},\bar{x})}(\bar{t}), v) \geq 0
	\]
hence the result. The last part of the proposition, in the differentiable case, may be found in \cite{Ester_Lev_arxiv}.
\end{proof}

\subsection{The Barrier Theorem}\label{barrierTh:subsec}

The barrier's construction is done according to the following: 
\begin{thm}\label{BarrierTheorem1}
	\sloppy Assume (A2.1) to (A2.5) hold. Consider an integral curve $x^{\bar{u}}$ on $\DAM \cap \cl(\Int(\AAA))$ and assume that the control function $\bar{u}$ is piecewise continuous. Then $\bar{u}$ and $x^{\bar{u}}$ satisfy the following necessary conditions.
	
	There exists a non-zero absolutely continuous adjoint $\lambda^{\bar{u}}$ and piecewise continuous multipliers $\mu_i^{\bar{u}} \geq 0$, $i=1,\dots,p$, such that:
	\begin{equation}\label{costateEquation1}
	\begin{aligned}
	\dot{\lambda}^{\bar{u}}(t) =-& \left(\frac{\partial f}{\partial x}(x^{\bar{u}}(t),\bar{u}(t)) \right)^T \lambda^{\bar{u}}(t) 
	\\&\quad \quad - \sum_{i=1}^{p}\mu_i^{\bar{u}}(t)\frac{\partial g_i}{\partial x}(x^{\bar{u}}(t),\bar{u}(t))
	\end{aligned}
	\end{equation}
	with the ``complementary slackness condition''
	\begin{equation}\label{eq:complement1}
	\mu_i^{\bar{u}}(t)g_i(x^{\bar{u}}(t),\bar{u}(t)) = 0, \quad i=1,\ldots, p.
	\end{equation}
	Moreover, at almost every $t$, the Hamiltonian, denoted by $\mathcal{H}(x^{\bar{u}}(t),u,\lambda^{\bar{u}}(t)) = \left(\lambda^{\bar{u}}(t) \right)^Tf(x^{\bar{u}}(t),u)$, is minimised over the set $U(x^{\bar{u}}(t))$ and constant:
	\begin{equation}\label{HamiltonianMinimised1}
	\begin{aligned}
		&\min_{u\in U(x^{\bar{u}}(t))} \lambda^{\bar{u}}(t)^T f(x^{\bar{u}}(t),u)\\ &=\min_{u\in U} \left[ \left(\lambda^{\bar{u}}(t) \right)^T f(x^{\bar{u}}(t),u) + \sum_{i=1}^{p} \mu_{i}^{\bar{u}}(t)g_{i}(x^{\bar{u}}(t),u)\right] \\
		&= \lambda^{\bar{u}}(t)^T f(x^{\bar{u}}(t),\bar{u}(t)) 
	\end{aligned}
	\end{equation}
	with the following boundary conditions:
	\begin{itemize}
	\item If the barrier ends on $G_0$ at a non differentiability point, then at this point the adjoint satisfies
		\begin{equation}\label{eq:finalConditions1}
		\lambda^{\bar{u}}(\bar{t}) = \left( D\tilde{g}(z_-) \right)^T
		\end{equation} 
		where $z$ and $\bar{t}$ are such that  $z= x^{\bar{u}}(\bar{t})\in G_0$ and \eqref{ineq:ult_tan}, namely
		$\min_{u\in U(z_-)} D\tilde{g}(z_-) f(z, u) \geq 0$,
		where $D\tilde{g}(z_-)$ indicates the left limit of $D\tilde{g}(x^{\bar{u}}(\tau))$, when $\tau \nearrow 			\bar{t}$.
	\item If the barrier ends on $G_0$ at a differentiability point, then at this point the adjoint satisfies
		\begin{equation}\label{eq:finalConditions2}
		\lambda^{\bar{u}}(\bar{t}) = \left( D\tilde{g}(z) \right)^T
		\end{equation} 
		where $z$ and $\bar{t}$ are such that  $z= x^{\bar{u}}(\bar{t})\in G_0$ and \eqref{eq:ult_tan_smooth}, namely
		$\min_{u\in U(z)} D\tilde{g}(z) f(z, u) = 0$
		\end{itemize}
\end{thm}
\begin{proof} Easy adaptation of the proof of Theorem 5.1 of \cite{Ester_Lev_arxiv}. See also Theorem 3.4.1 of \cite{EsterhuizenPHDThesis}.
\end{proof}

The reader may find a thorough discussion of this result, its limitations and related open problems in \cite{EsterhuizenPHDThesis}.

\section{The Barrier for the Pendulum on a Cart with a Non-Rigid Cable}\label{sec:Barrier_for_Pend}
\subsection{Derivation of Constrained System}\label{subsec:Derivation_Of_Const_Sys}

\begin{figure}
	\begin{center}
		\includegraphics[height=4cm]{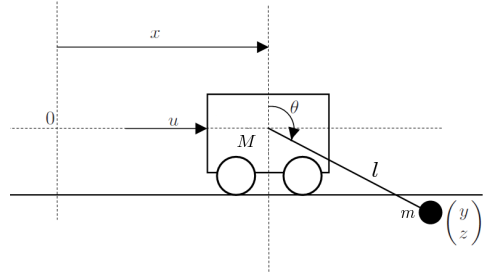}
		\caption{Pendulum on a cart with non-rigid cable.}
		\label{Fig:PendulumOnCart}
	\end{center}
\end{figure}

We consider the system as in Figure \ref{Fig:PendulumOnCart}: a mass of $m$ (kg) is attached to the end of a massless cable that may go slack, which is suspended from a cart of mass $M$ (kg) that may move unconstrained along a horizontal line. The control, $u$, is the force (N) applied to the cart satisfying $|u| \leq 1$. The angle (rd) between the cable and the vertical is $\theta$, $l$ is the length (m) of the cable and $g$ the acceleration due to gravity. The cart's position is $x$, and the coordinates of the mass $m$ are $(y,z)$. As long as the cable is taut, $l$ is constant, $y = x + l\sin \theta$ and $z = l\cos \theta$.

One way to guarantee tautness in the cable is to impose the condition that the cable's tension, $T$, is always nonnegative. Under this assumption the dynamics of the system, obtained via the Euler-Lagrange method, are given by:
\begin{align}
\label{eq:pend1}
& \dot{\theta}_1  =  \theta_2 , \\
\label{eq:pend2}
& \dot{\theta}_2 = \frac{-u\cos\theta_1 + (M + m)g\sin\theta_1 - ml\theta^2_2\cos\theta_1\sin\theta_1}{l\left(M + m\sin^2\theta_1\right)}  \\
\label{eq:pend3}
&\dot{x}_1 = x_2 \\
\label{eq:pend4}
&\dot{x}_2 = \frac{u + ml\theta^2_2\sin\theta_1 - mg\cos\theta_1\sin\theta_1}{M + m\sin^2\theta_1} 
\end{align}
where $x_1 = x$, and $\theta_1 = \theta$. To lighten our notation, we introduce $\qq \triangleq (\theta_1,\theta_2,x_1,x_2)$ and $\tth \triangleq (\theta_1,\theta_2)$.
Remark that the dynamics \eqref{eq:pend1}, \eqref{eq:pend2} of $\tth$, where no simplification or approximation of any kind has been made, do not  depend explicitly on the cart's position and velocity $(x_1,x_2)$ and that the dynamics of $(x_1,x_2)$ and $\tth$ are only coupled via the force $u$.

We now show that imposing the condition that the tension in the cable remains nonnegative is equivalent to imposing a \emph{mixed} constraint on the system.

By considering the balance of forces on the mass $m$, its projection on the vertical axis is indeed given by $m\ddot{z}= -mg -T\cos\theta_1$ (see for example \cite{JLbook}). Thus $T= - \frac{m \left(\ddot{z}+g\right)}{\cos\theta_{1}} \geq 0$, which is equivalent to:
\begin{equation}\label{eq:balance_of_forces}
-\frac{\ddot{z} + g}{\cos \theta_1} \geq 0.
\end{equation}
Noting that $\ddot{z} = -l\left(\theta_2^2 \cos \theta_1 + \dot{\theta}_2 \sin \theta_1 \right)$, we substitute equation \eqref{eq:pend2} in the latter expression and multiply \eqref{eq:balance_of_forces} by $l\left(M + m\sin^2 \theta_1 \right)$. The inequality then simplifies to the mixed constraint:
\begin{equation}\label{ineq:PendulumMixedConstraint}
u\sin\theta_1 + Mg\cos\theta_1 - Ml\theta_2^2 \leq 0.
\end{equation}
Thus, the problem is to obtain the barrier for the system described by \eqref{eq:pend1}-\eqref{eq:pend4} subject to the constraint on the control, $|u| \leq 1$, and the mixed constraint \eqref{ineq:PendulumMixedConstraint}.  We do not consider a constraint on the cart's track length for clarity's sake.
Note that zero tension in the cable results in free-fall but, depending on the cart's trajectory, the cable may remain taut or become slack. 

Note that the smoothness and convexity assumptions (A2.1), (A2.3), (A2.4), are indeed satisfied by the right-hand side of \eqref{eq:pend1}--\eqref{eq:pend4} and constraint \eqref{ineq:PendulumMixedConstraint} and that (A2.2) is only satisfied for bounded $\theta_2$, which is not a real limitation in vue of Remark~\ref{rem:A22}.

\subsection{Constructing the Barrier}

\begin{figure}[thpb]
	\begin{center}
		\includegraphics[height=12cm]{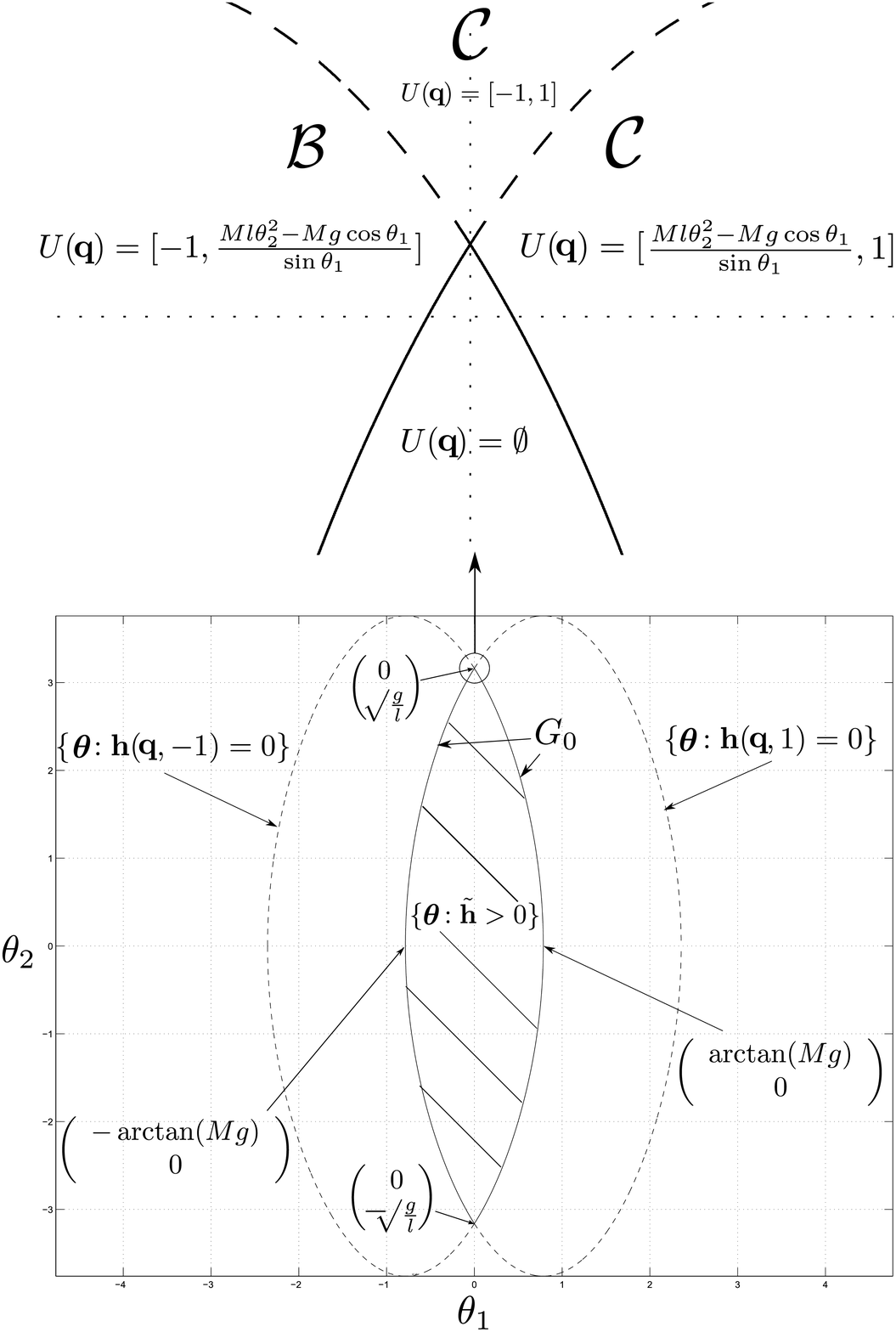}
		\caption{The set $\{\boldsymbol{\theta} : -|\sin\theta_1| + Mg\cos\theta_1 - Ml\theta_2^2 = 0 \}$ for $M = 0.1$, $m = 0.1$, $l = 1$ and $g = 10$, along with some important points. The top figure presents a closer look at the point $(0,\sqrt{\frac{g}{l}})$ and specifies the set $U(\textbf{q})$ in various parts of a neighbourhood of the point.}
		\label{Fig:HTilde_Pendulum}
	\end{center}
\end{figure}

Recalling the notations $\qq = (\theta_1,\theta_2,x_1,x_2)$ and $\tth = (\theta_1,\theta_2)$, we label the mixed constraint $\hh (\textbf{q},u) =u\sin\theta_1 + Mg\cos\theta_1 - Ml\theta_2^2$ and so 
$$\tilde{\hh}(\textbf{q}) = \min_{|u|\leq 1} \mathbf{h}(\textbf{q},u) = -|\sin\theta_1| + Mg\cos\theta_1 - Ml\theta_2^2.$$ Let us assume that barrier trajectories reach the set $G_0= \{\qq : \tilde{\hh}(\qq) = 0  \}$, whose projection onto the plane $(\theta_1,\theta_2)$ is shown in Figure \ref{Fig:HTilde_Pendulum}. 

Note that the equation $\tilde{\hh}(\qq) = 0$ only has a solution for $\theta_1\in[-\arctan(Mg) + 2k\pi,\arctan(Mg) + 2k\pi]$, $k\in \NN$, and that $\tilde{\mathbf{h}}$ is not differentiable if $\theta_1 = 2k\pi$. Also, according to \eqref{def:U(x)}, we have
\begin{equation}\label{eq:U(x)}
U(\textbf{q})
=\begin{cases}
[-1, \frac{Ml\theta_2^2 - Mg\cos\theta_1}{\sin\theta_1}]\cap[-1,1] & \mathrm{if}\,\, \sin\theta_1 > 0\\
[\frac{Ml\theta_2^2 - Mg\cos\theta_1}{\sin\theta_1}, 1]\cap[-1,1] & \mathrm{if}\,\, \sin\theta_1 < 0\\
[-1,1] & \mathrm{if}\,\, \sin\theta_1 = 0.
\end{cases}
\end{equation}
It is easily verified that the independence condition of assumption (A2.5) is met everywhere in a neighborhood of $G_0$, except on $G_0$ itself, where $U(\qq)=\{\pm1\}$ and $\#\II(\qq,u)= 2$. Therefore (A2.5) is satisfied.

\subsubsection{Barrier End Points}\label{subsec:Pts_Of_Tan}

\sloppy We first look at points on $G_0$ where $\tilde{\hh}$ is differentiable and without loss of generality we will only carry out the analysis for $\theta_1\in[-\arctan(Mg), 0[$. 

Invoking equation \eqref{eq:ult_tan_smooth} as well as the final condition \eqref{eq:finalConditions2}, we obtain:
\begin{equation}\label{eq:Ult_Tan_pend}
\begin{array}{l}
\min_{u\in\{ 1 \}} \left(\cos\theta_1 - Mg\sin\theta_1\right)\theta_2 \\
+ 2Ml\theta_2\left(\frac{u\cos\theta_1 - (M+m)gsin\theta_1 + ml\theta_2^2\cos\theta_1\sin\theta_1}{l\left(M + m\sin^2(\theta_1)\right)}\right) = 0
\end{array}
\end{equation}
where $U(\textbf{q}) = \{1 \}$ since, for $\textbf{q}\in G_0$,   $\lim_{(\theta_{1}, \theta_{2}) \rightarrow (-\arctan Mg,0)}\frac{Ml\theta_2^2 - Mg\cos\theta_1}{\sin\theta_1} = 1$. From here we easily identify $(\theta_1^{\bar{u}}(\bar{t}),\theta_2^{\bar{u}}(\bar{t})) = (-\arctan(Mg),0)$, with $x_1^{\bar{u}}(\bar{t})$ and $x_2^{\bar{u}}(\bar{t})$ free, as end points, along with the final adjoint given by \eqref{eq:ult_tan_smooth}:
\begin{equation}
\begin{array}{l}
\ds \lambda^{\bar{u}}(\bar{t})  = \\ \hspace{0.5cm}\ds (\cos\theta_1^{\bar{u}}(\bar{t}) - Mg\sin\theta_1^{\bar{u}}(\bar{t}),-2Ml\theta_2^{\bar{u}}(\bar{t}),0,0)^T.\label{eq:final_costate_smooth}
\end{array}
\end{equation}
Let us show that \eqref{eq:Ult_Tan_pend} does not have another solution for any $\theta_1\in[-\arctan(Mg),0[$. Indeed, $\theta_2 \neq 0$ and we must investigate:
\begin{equation}
\begin{array}{l}
\left(\cos\theta_1 - Mg\sin\theta_1\right) \\
+ 2Ml\left(\frac{\cos\theta_1 - (M+m)gsin\theta_1 + ml\theta_2^2\cos\theta_1\sin\theta_1}{l\left(M + m\sin^2(\theta_1)\right)}\right) = 0.
\end{array}
\end{equation}
We now substitute $\theta_2^2$ using $\tilde{\mathbf{h}}(\textbf{q}) = 0$, multiply by $l\left(M + m\sin^2(\theta_1)\right)$ and use the identity $m\cos^3\theta_1 - m\cos\theta_1 = -m\cos\theta_1\sin^2\theta_1$ to arrive, after some algebra, at the expression:
\begin{equation}
\begin{array}{l}
-M\cos\theta_1 - Mmg\sin\theta_1\cos^2\theta_1\\
+ M^2g\sin\theta_1 + Mmg\sin\theta_1 - \sin^2\theta_1(m\cos\theta_1) = 0.
\end{array}
\end{equation}
After grouping terms we get:
\[
\left( Mg\sin\theta_1 -\cos\theta_1 \right) \left(M + m\sin^2\theta_1\right) = 0.
\]
Since $\left(M + m\sin^2\theta_1\right)>0$ we get $\theta_1 = \arctan(\frac{1}{Mg}) \notin[-\arctan(Mg),0[$, and so there is not another solution for $\theta_1\in [-\arctan(Mg),0[$.

Along the same lines we deduce that all the points $(\theta_1,\theta_2) = (\pm\arctan(Mg) + 2k\pi,0)$, $k\in \NN$, with $x_1$ and $x_2$ free, are the only end points on $G_0$ where $\tilde{\hh}$ is differentiable.

We now turn our attention to the point $(\theta_1,\theta_2) = (0,\sqrt{\frac{g}{l}})$ (with $x_1$ and $x_2$ arbitrary) where $\tilde{\mathbf{h}}$ is not differentiable; the analysis will carry over to the points $(\theta_1,\theta_2) = (2k\pi,\pm\sqrt{\frac{g}{l}})$, $k\in \NN$, in a similar way. We introduce the following sets:
\[
\mathcal{B} \triangleq \left\{\boldsymbol{\theta} : \theta_1 < 0, \frac{Ml\theta_2^2 - Mg\cos\theta_1}{\sin\theta_1} \geq -1, \tilde{\mathbf{h}}(\textbf{q}) \leq 0 \right\}
\]
\[
\mathcal{C} = \mathcal{B}^{\mathsf{C}} \setminus \{\boldsymbol{\theta} : \tilde{\mathbf{h}}(\textbf{q}) > 0 \}
\]
as in Figure \ref{Fig:HTilde_Pendulum}. From Theorem \ref{BarrierTheorem1}, if a barrier trajectory intersects the point $(\theta_1,\theta_2) = (0,\sqrt{\frac{g}{l}})$, with $x_1$ and $x_2$ arbitrary, at time $\bar{t}$ (without confusion we use the same label for this time instant as was used previously for the analysis of the point $(\theta_1,\theta_2) = (-\arctan(Mg),0)$) then condition \eqref{ineq:ult_tan} holds. If this barrier trajectory approached the point from the set $\mathcal{C}$ then it can be verified that we would get
\[
\min_{u\in U(\textbf{q}^{\bar{u}}(\bar{t}_-))} D\tilde{\mathbf{h}}(\textbf{q}^{\bar{u}}(\bar{t}_-)) f(\textbf{q}^{\bar{u}}(\bar{t}), u)< 0,
\]
where $D\tilde{\hh}$  denotes the gradient of $\tilde{\hh}$ with respect to the vector $\qq$ and with $f$ the right-hand side of \eqref{eq:pend1}--\eqref{eq:pend4},
which would violate condition \eqref{ineq:ult_tan}. Moreover, this trajectory can clearly not approach the point from the set $\{\theta : \tilde{\hh}(\qq) > 0 \}$. The only possibility left is that it approaches the point from the set labelled $\mathcal{B}$, and the final adjoint is given by \eqref{eq:finalConditions1} with condition \eqref{ineq:ult_tan}:
\begin{equation}\label{eq:costate_non_diff}
\lambda^{\bar{u}}(\bar{t}_-)^T = D\tilde{\hh}(\qq^{\bar{u}}(\bar{t}_-)) = (1, -2Ml\sqrt\frac{g}{l},0,0).
\end{equation}

\subsubsection{Deriving the Control Associated with the Barrier}

The adjoint equations are given by \eqref{costateEquation1}, from which it can be verified that $\dot{\lambda}_3 = 0$ and $\dot{\lambda}_4 = \lambda_3$.

From equations \eqref{eq:final_costate_smooth} and \eqref{eq:costate_non_diff}, we deduce that $\lambda_3(t) = \lambda_4(t) \equiv 0$. The Hamiltonian is here given by $$\begin{aligned} &\mathcal{H}(\textbf{q},u,\lambda) = \lambda_1 \theta_2 \\ &\quad + \lambda_2\frac{-u\cos\theta_1 + (M + m)g\sin\theta_1 - ml\theta^2_2\cos\theta_1\sin\theta_1}{l\left(M + m\sin^2\theta_1\right)}\end{aligned}$$ and from
\eqref{HamiltonianMinimised1} we have that $\mu$ satisfies: $\frac{\partial H}{\partial u} + \mu\frac{\partial \mathbf{h}}{\partial u} = 0$
which gives:
\[
\mu(t) = \left\{\begin{array}{ll} \frac{\lambda_2(t)\cot\theta_1(t)}{l\left(M + m\sin^2(\theta_1)\right)} & ~\mathrm{if}~\mathbf{h}(\textbf{q}(t),\bar{u}(t)) = 0\\
0 & ~\mathrm{if}~ \mathbf{h}(\textbf{q}(t),\bar{u}(t)) < 0.
\end{array}\right.
\]
The Hamiltonian mimimisation yields:
\[
\begin{array}{lll}
\mathrm{if}~\lambda_2(t)\cos\theta_1(t) > 0 & \\
\,\,\,\,\bar{u}(t)\hspace{-2pt}=\hspace{-3pt}\left\{
\begin{array}{lll}
\hspace{-2pt}1 & \mathrm{if} & \sin\theta_1(t)\leq 0\\
\hspace{-2pt}\min \left( \frac{Ml\theta_2^2(t) - Mg\cos\theta_1(t)}{\sin\theta_1(t)} , 1\right)& \mathrm{if} & \sin\theta_1(t) > 0
\end{array}\right. &\\
\mathrm{if}~\lambda_2(t)\cos\theta_1(t) < 0 & \\
\,\,\,\,\bar{u}(t)\hspace{-2pt}=\hspace{-3pt}\left\{
\begin{array}{lll}
\hspace{-2pt}-1 & \hspace{-2pt}\mathrm{if} & \sin\theta_1(t) \geq 0\\
\hspace{-2pt}\max\left(\frac{Ml\theta_2^2(t) - Mg\cos\theta_1(t)}{\sin\theta_1(t)},-1 \right)& \hspace{-2pt}\mathrm{if} & \sin\theta_1(t) < 0
\end{array}\right. &\\
\mathrm{if}~\lambda_2(t)\cos\theta_1(t) = 0 & \\
\,\,\,\,\bar{u}(t)=\mathrm{arbitrary}.&
\end{array}
\]

\subsubsection{Backward Integration of System Equations}\label{subsection:Backward_Integration}

As previously remarked, the right hand sides of equations \eqref{eq:pend1} and \eqref{eq:pend2}, as well as the mixed constraint, equation \eqref{ineq:PendulumMixedConstraint}, and the control associated with the barrier, $\bar{u}$, are independent of $x_1$ and $x_2$. Moreover, as shown in subsection \ref{subsec:Pts_Of_Tan}, the values of $x_1$ and $x_2$ are arbitrary at points of ultimate tangentiality. These facts allow us to simplify the analysis by ignoring $x_1$ and $x_2$ from this point forward, only focusing on $\tth =(\theta_1,\theta_2)$.

\begin{figure}[thpb]
	\begin{center}
		\includegraphics[height=7.9cm]{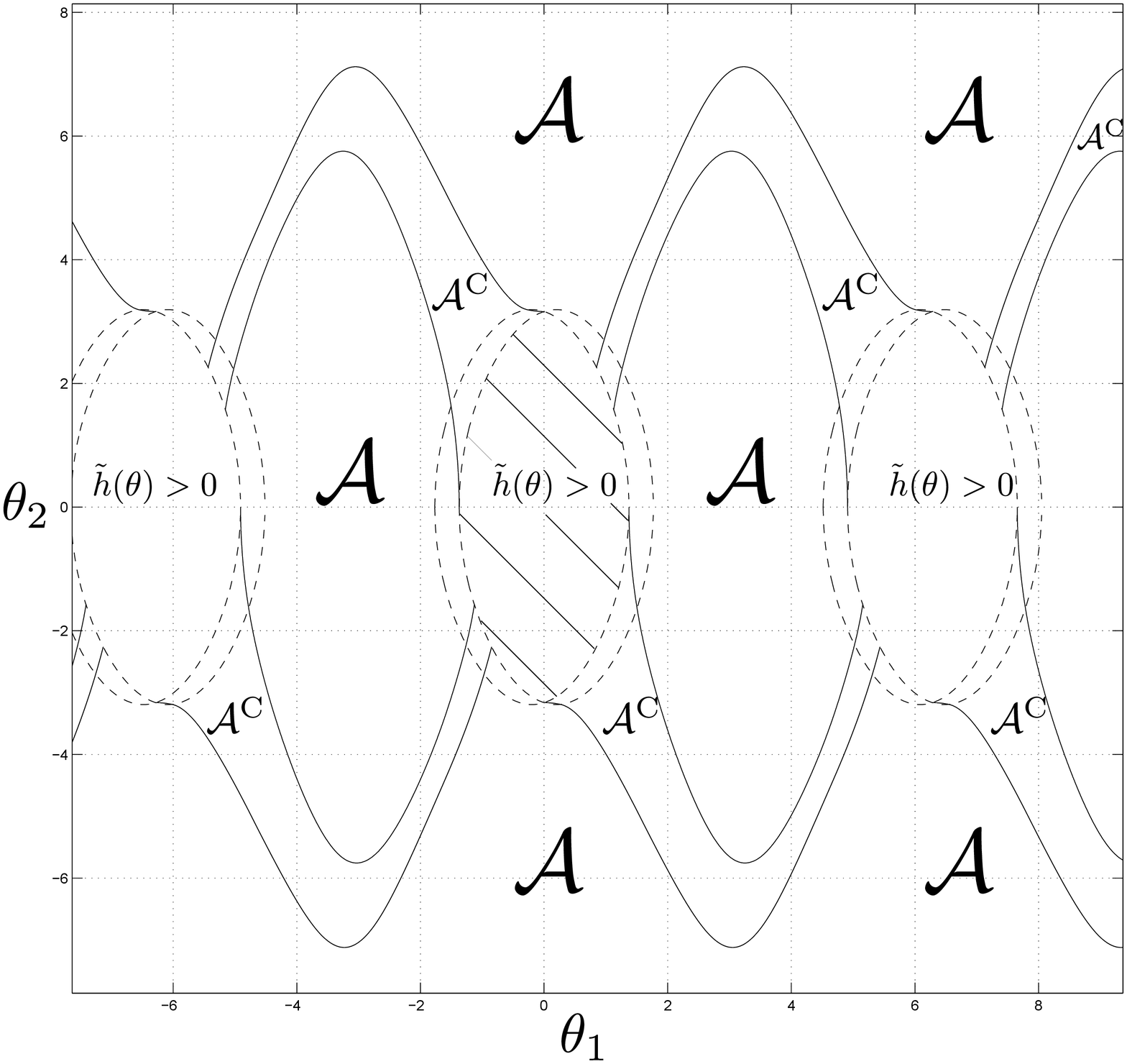}
		\caption{The admissible set for the pendulum on a cart with a nonrigid cable, with the constraint that the tension in the cable is always nonnegative, \eqref{ineq:PendulumMixedConstraint}. The constants in this case are: $g = 10\mathrm{m/s}^2$, $l = 1$ metre, $M = 0.5$ kg, $m = 0.1$ kg. Note that the admissible set is made up of disjoint parts.} 
		\label{Fig:BigMassPendulum}
	\end{center}
\end{figure}

\begin{figure}[thpb]
	\begin{center}
		\includegraphics[width=8.3cm]{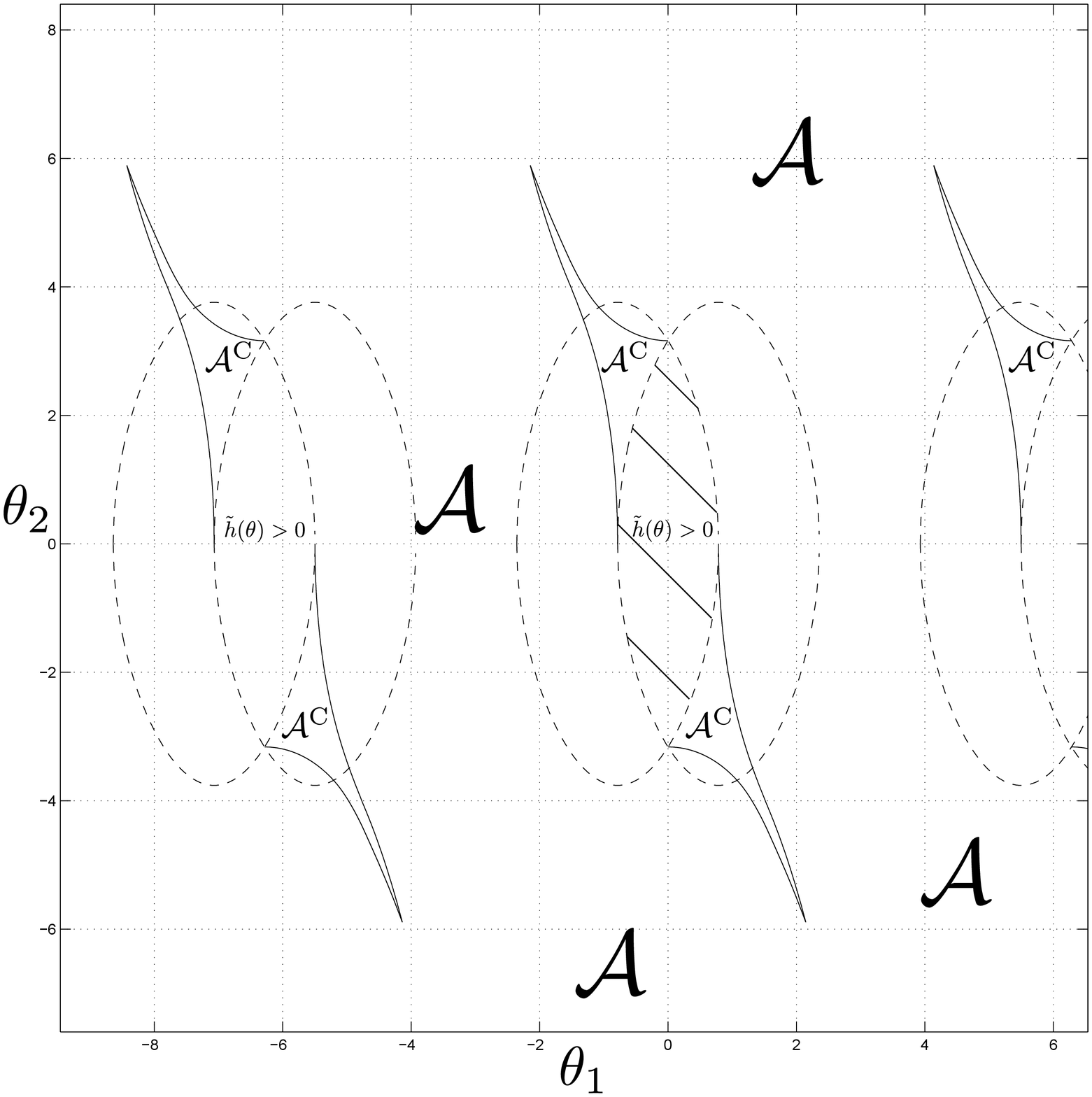}
		\caption{The admissible set for the pendulum on a cart with a slack rope, equations, with the constraint that the tension in the cable is always nonnegative, \eqref{ineq:PendulumMixedConstraint}. The constants in this case are: $g = 10\mathrm{m/s}^2$, $l = 1$ metre, $M = 0.1$ kg, $m = 0.1$ kg. Note that the obtained barrier trajectories intersect at \emph{stopping points}.} 
		\label{Fig:SmallCartMass}
	\end{center}
\end{figure}

Integrating equations \eqref{eq:pend1}, \eqref{eq:pend2} and \eqref{costateEquation1} backwards using $\bar{u}$, the expression of \eqref{costateEquation1} being omitted for clarity's sake, from the identified points of ultimate tangentiality, $(\pm\arctan(Mg) + 2k\pi,0)$ and $(2k\pi,\pm\sqrt{\frac{g}{l}})$, $k\in \NN$, utilising the appropriate final adjoint for each point, we obtain the trajectories as in Figure \ref{Fig:BigMassPendulum} and Figure \ref{Fig:SmallCartMass}.

\begin{figure}[thpb]
	\begin{center}
		\includegraphics[width=0.8\columnwidth]{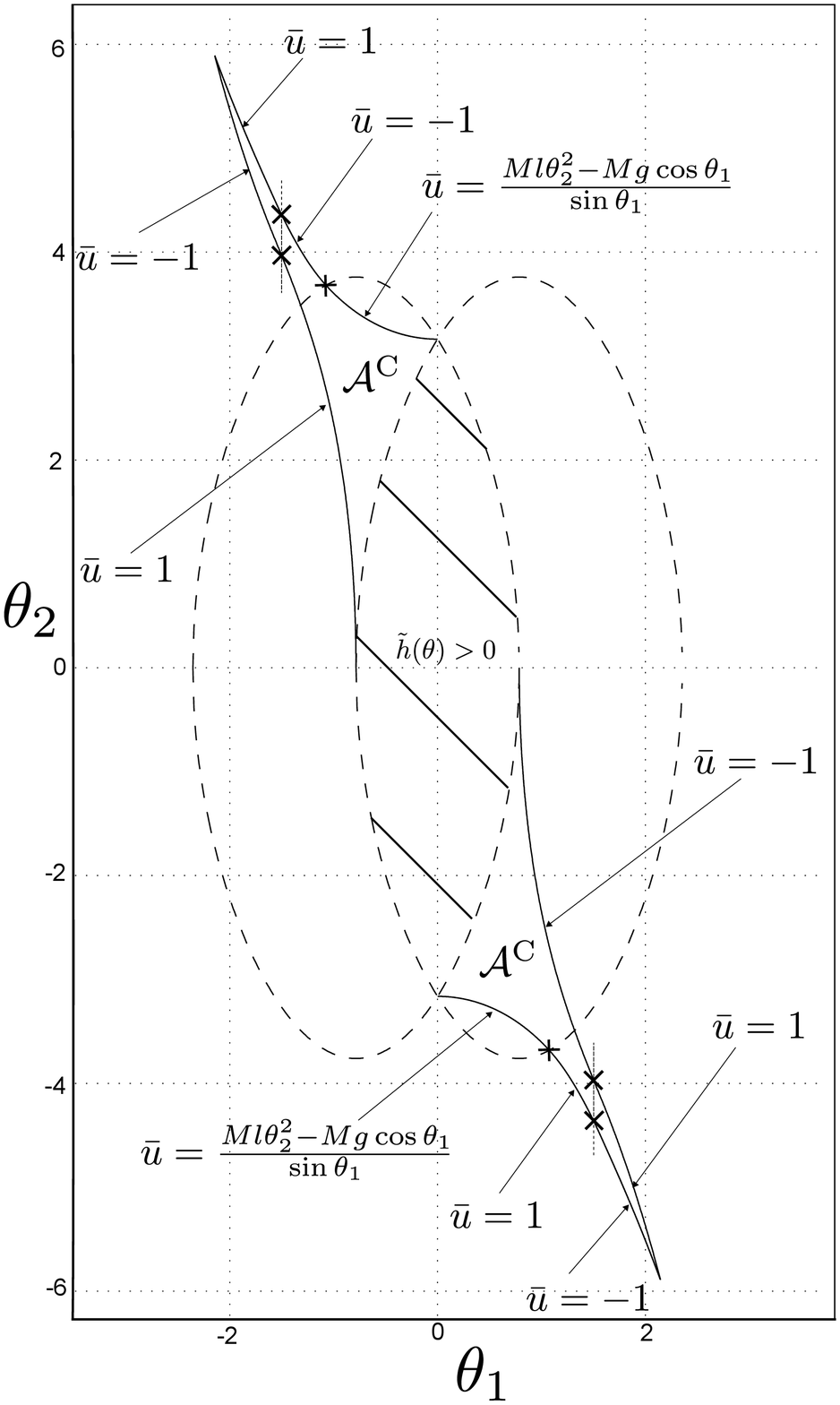}
		\caption{A closer look at the control associated with the barrier trajectories from Figure \ref{Fig:SmallCartMass}. The ``$\times$'' correspond to $\theta_1 = \pm\frac{\pi}{2}$ where the control switches. The ``+'' correspond to the points where the controls associated with the barrier trajectories arriving at $(0,\pm\frac{g}{l})$ switch to $\frac{Ml\theta_2^2 - Mg\cos\theta_1}{\sin\theta_1}$.}
		\label{Fig:SmallCartMassZoomedIn}
	\end{center}
\end{figure}

Figure \ref{Fig:BigMassPendulum} shows the obtained admissible set for a set of constants where the mass of the cart is big relative to the pendulum mass. On the contrary, when the cart mass is relatively small, the barrier trajectories intersect and by Theorem \ref{thm:stopping_points} of Appendix \ref{Appendix:stop}, we deduce that these intersection points are stopping points, see Figure \ref{Fig:SmallCartMass}. Figure \ref{Fig:SmallCartMassZoomedIn} provides a closer look at the control function along barrier trajectories as obtained for the constants in Figure~\ref{Fig:SmallCartMass}.

\section{Discussion}\label{sec:Discussion}
If $\boldsymbol{\theta}$ is such that $\tilde{\mathbf{h}}(\textbf{q}) \geq 0$ then the angular velocity of the mass is too small to provide positive tension in the cable, and the mass enters free-fall. If $\boldsymbol{\theta}\in \AC$ then no admissible control function can prevent a loss of tautness in the future.

At the points of ultimate tangentiality $(\pm\arctan(Mg) + 2k\pi, 0)$, $k\in \NN$, the angular velocity of the mass is zero as well as the tension in the cable and the mass is in free-fall. However, employing the only admissible control at this point (i.e. $u = \pm 1$ depending on the point) results in the state immediately entering the interior of the admissible set and the tension can be made positive again.

At the singular points $(2k\pi,\pm \sqrt\frac{g}{l})$, $k\in \NN$, the control acts perpendicularly to the cable and so does not have any effect on the tension. The barrier trajectory that passes through these points is quite interesting because along the entire part of the trajectory for which $\bar{u} = \frac{Ml\theta_2^2 - Mg\cos\theta_1}{\sin\theta_1}$ the mass is in free-fall but taut. When the barrier arrives at the points $(2k\pi,\pm \sqrt\frac{g}{l})$ it is again possible to employ a control such the state enters the interior of the admissible set and the tension becomes positive again.

If the constants are such as those specified in Figure \ref{Fig:BigMassPendulum}  the admissible set consists of a periodic sequence of two connected components, one of them being bounded and the other not. This has the interpretation that if the system is initiated in the bounded one then it is impossible to increase the angular velocity beyond a certain bound without leaving the admissible set. On the contrary, if the system is initiated in the unbounded one, one can ``spin'' the mass through the full range of angles, i.e. through all $\theta_1 \in \RR$, always maintaining a taut cable. Let us stress that no control allows the state to pass from one component to the other without entering $\AC$.

The admissible set obtained in Figure \ref{Fig:SmallCartMass} is now connected thus allowing more manoeuvrability.  

\section{Concluding Remarks and Future Research}\label{sec:Conclusions}

We can model the pendulum on a cart with a non-rigid cable as a hybrid automaton, see for e.g. \cite{VanDerShaft2000}. In this framework a hybrid system is specified by a graph with the nodes corresponding to ``locations'', where at each location the continuous state evolves according to a particular differential equation. At ``event times'' there occur ``transitions'' between locations along the graph's edges.

The pendulum system may be modelled as a hybrid automaton with two locations: at the first location the continuous state evolves according to \eqref{eq:pend1} - \eqref{eq:pend4} and at the second location the state evolves according to the free-fall dynamics of the mass with slack cable. The admissible set can then be interpreted as a \emph{potentially safe set}, in other words if the state remains in this set it is guaranteed that there exists a control such that the system does not transition out of the initial location. We contrast this with the usual notion of \emph{safety sets} in hybrid systems where it is generally required that the state remains in this set \emph{for all} possible control functions, see e.g. \cite{Mitchell2003} and \cite{Lygeros2007}, hence the term \emph{potentially} safe.

The results may find application in the study of other similar problems in engineering, such as the control of weight-handling equipment, UAVs and to obtaining potentially safe sets in hybrid systems. Indeed, the same approach should be applicable to higher dimensional problems such as the pendulum in 3 dimensions with non-rigid cable. This application will be the subject of future research.

Future research could also focus on the development of a richer theory of potentially safe sets for hybrid automata with any finite number of locations, similar to the ideas in \cite{Mitchell2003}.

\bibliographystyle{plain}        
\bibliography{autosam1}          

\appendix
\section{Barrier Stopping Points}\label{Appendix:stop}    

Backwards integrated barrier trajectories obtained from Theorem \ref{BarrierTheorem1} may intersect, with their further backward prolongations being in the interior of the admissible set. In this case these prolongations need to be ignored. A preliminary study of this phenomenon has been presented in \cite{EL_IFACE2014} and we summarise the main result that we  use in the construction of the barrier for the pendulum problem.
\begin{definition}\label{def:StoppingPointDef}
	Consider two distinct integral curves $x^{(u_1,z_1)}$ and $x^{(u_2,z_2)}$ obtained from Theorem \ref{BarrierTheorem1} by backward integration, running along the barrier $\DAM$ from two distinct points $z_1,z_2\in G_0$ at $\bar{t}_1$ and $\bar{t}_2$ respectively, i.e. $x^{(u_i,z_i)}(\bar{t}_i) = z_i$, $i = 1,2$, where $u_i$ is the corresponding control function that satisfies condition \eqref{HamiltonianMinimised1} for almost all $t\leq \bar{t}_i$, $i = 1,2$. Assume that there exists a point of transversal \footnote{in other words with $f(\xi,u_1(\tilde{t}))$ and $f(\xi,u_2(\tilde{t}))$ independent} intersection $\xi$ of these two curves at some time labeled $\tilde{t}$. $\xi$ is said to be a \emph{barrier stopping point by intersection} either if the two maximal integral curves stop at $\xi$, or if $x^{(u_i,z_i)}(t)\in\Int(\AAA)$, $i = 1,2$, for all $t<\tilde{t}$, whereas $x^{(u_i,z_i)}(t)\in\DAM$ for all $t\in[\tilde{t},\bar{t}_i]$, $i = 1,2$.
\end{definition}
\begin{thm}\label{thm:stopping_points}
	Consider two distinct integral curves $x^{(u_1,z_1)}$ and $x^{(u_2,z_2)}$ as in Definition \ref{def:StoppingPointDef}. If there exists an intersection point $\xi$ of these two curves at some time\footnote{in case of multiple intersection points, only the largest time $\tilde{t}<\bar{t}_i$, $i = 1,2$, must be considered.} $\tilde{t}$, i.e. $x^{(u_1,z_1)}(\tilde{t}) = x^{(u_2,z_2)}(\tilde{t}) = \xi$, then $\xi$ is a \emph{barrier stopping point by intersection}.
\end{thm}

\section{Needle Perturbations and the Variational Equation}\label{Appendix:var}

Given $\bar{u} \in \UU$ and an integral curve $x^{(\bar{u},\bar{x})}$, we consider $\ee > 0$ and bounded, an initial state perturbation $h\in \RR^n$ satisfying $\Vert h\Vert \leq H$ and a variation $u_{\kappa,\ee}$ of $\bar{u}$, parameterised by the vector $\kappa \triangleq (v,\tau,l) \in U(x^{(\bar{u},\bar{x} + \ee h)}(\tau - \ee l)) \times [0,T] \times [0,L]$ with bounded $T,L$, of the form
\begin{equation}\label{u-var-eq}
u_{\kappa,\ee} \triangleq 
\bar{u} \Join_{(\tau-l\ee)} v \Join_{\tau} \bar{u}
=
\left\{ \begin{array}{ll}
v&\mbox{\textrm{on}~} [\tau-l\ee, \tau[\\
\bar{u}&\mbox{\textrm{elsewhere on}~} [0,T]
\end{array}\right. 
\end{equation}
where $v$ stands for the constant control equal to $v \in U(x^{(\bar{u},x_0)}(\tau))$ for all $t\in [\tau-l\ee, \tau[$. 
We have $x^{(u_{\kappa,\ee}, \bar{x}+\ee h)}(t)= x^{(\bar{u},\bar{x}+\ee h)}(t)$ for all $t\in [0, \tau-l\ee [$ and, denoting by $z_{\ee} (\tau-l\ee) \triangleq  x^{(\bar{u},\bar{x}+\ee h)}(\tau-l\ee)$ and 
$z_{\ee}(\tau) \triangleq x^{(u_{\kappa,\ee}, \bar{x}+\ee h)}(\tau)$, we have 
\begin{equation}
	\begin{aligned}
&x^{(u_{\kappa,\ee}, \bar{x}+\ee h)} =\\
& x^{(\bar{u},\bar{x}+\ee h)} \Join_{(\tau-l\ee)} x^{(v, z_{\ee} (\tau-l\ee), \tau - l\ee)} \Join_{\tau} x^{(\bar{u}, z_{\ee}(\tau),\tau)}
	\end{aligned}
\end{equation}
We also consider the fundamental matrix of the variational equation ($I_{n}$ is the identity matrix of $\RR^{n}$):
\begin{equation}\label{group}
	\begin{aligned}
		&\frac{d}{dt}\Phi^{\bar{u}}(t,s) = \left( \frac{\partial f}{\partial x}(x^{(\bar{u},\bar{x})}(t),\bar{u}(t))\right) \Phi^{\bar{u}}(t,s),\\ & \Phi^{\bar{u}}(s,s)=I_{n}.
	\end{aligned}
\end{equation}
\begin{lem}\label{approx-lem} \emph{\cite[Chapter II, \S 13]{PBGM}, \cite[Chapter 4, p. 248]{Lee_Markus}, \cite{Gam}]} 
	The sequence $\left\{ x^{(u_{\kappa,\ee}, \bar{x}+\ee h)} \right\}_{\ee \geq 0}$ is uniformly convergent  to $x^{(\bar{u},\bar{x})}$ on $[0,T]$ as $\ee\rightarrow 0$, uniformly with respect to $\kappa$ and $h$.
If, moreover,  $\tau$ is a Lebesgue point of $\bar{u}$, we have, for all $t\in [\tau, T]$:
	\begin{equation}\label{approx-eq}
	x^{(u_{\kappa,\ee}, \bar{x}+\ee h)}(t) - x^{(\bar{u},\bar{x})}(t) = 
	\ee w(t,\kappa,h)
	+ O(\ee^{2})
	\end{equation} 
\end{lem}
where
\begin{equation}\label{needle-eq}
	\begin{aligned}
	&w(t,\kappa,h) \triangleq \Phi^{\bar{u}}(t,0)h \\ & \,\,+ l \Phi^{\bar{u}}(t,\tau) \left( f(x^{(\bar{u},\bar{x})}(\tau), v) - f(x^{(\bar{u},\bar{x})}(\tau), \bar{u}(\tau)) \right).
	\end{aligned}
\end{equation}
\end{document}